\documentclass[]{article}
\usepackage{amsmath}
\usepackage{amssymb}
\usepackage{hyperref}
\usepackage{pgf,pgfarrows,pgfnodes,pgfautomata,pgfheaps,pgfshade}
\usepackage{graphicx,amsfonts}
\usepackage{amsmath,amssymb,setspace}
\usepackage{mathrsfs}
\usepackage{colortbl}
\usepackage{tikz}
\usepackage{framed}
\usepackage{subfig}

\numberwithin{equation}{section}

\newtheorem{thm}{Theorem}[section]

\newtheorem{lemma}[thm]{Lemma}

\newtheorem{remark}[thm]{Remark}

\newtheorem{definition}[thm]{Definition}

\newcommand{\R}{\mathbb{R}}

\title{Iterative Methods for Symmetric Outer Product Tensor Decompositions}

\author{Na Li\footnote{Department of Mathematics, Clarkson University, Potsdam, NY 13699, USA, nali@clarkson.edu.}  \  and Carmeliza Navasca\footnote{Department of Mathematics, University of Alabama at Birmingham, 1300 University Boulevard, Birmingham, AL, 35294, USA, cnavasca@uab.edu}}

\date{\today}

\begin{document}

\maketitle

\begin{abstract}
\setcounter{section}{0}

We study the symmetric outer product decomposition which decomposes a fully (partially) symmetric tensor into a sum of rank-one fully (partially) symmetric tensors. We present iterative algorithms for the third-order partially symmetric tensor and fourth-order fully symmetric tensor. The numerical examples indicate a faster convergence rate for the new algorithms than the standard method of alternating least squares.


\vspace{8pt}
\noindent
\end{abstract}

\section{Introduction}

In 1927, Hitchcock \cite{Hitch1}\cite{Hitch2} proposed the idea of the polyadic form of a tensor, i.e., expressing a tensor as the sum of a finite number of rank-one tensors. Today, this decomposition is called the canonical polyadic (CP); it is known as CANDECOMP or PARAFAC. It has been extensively applied to many problems in various engineering \cite{Sid1,Sid2,Acar,DeVos} and science \cite{Smilde,Kroonenberg}. Symmetric tensors have been used in many signal processing applications \cite{Comon1,Comon3,Lieven}. Similar with the CP decomposition for a general tensor, the symmetric outer product decomposition (SOPD) for fully symmetric tensors factors a fully symmetric tensor into a number of rank-one fully symmetric tensors. It related to the independent component analysis (ICA) \cite{Aapo,Comon2} or blind source separation (BSS), which is used to separate the true signal from noise and interferences in signal processing \cite{Comon3,Lieven}. For the SOPD of partially symmetric tensors, when the tensor order is 3 and it is symmetric on mode one and mode two, such a problem corresponds to the Indscal model introduced by Carrol and Chang \cite{JJ,AStegeman}. 

The well-known iterative method for implementing the sum of rank one terms is the Alternating Least-Squares (ALS) technique. Independently, the ALS was introduced by Carrol and Chang \cite{JJ} and Harshman \cite{RAH} in 1970. Since the SOPD is a special case of CP decomposition, the ALS method can be applied to solve the SOPD. A different method proposed by Comon \cite{Comon4} for SOPD reduces the problem to the decomposition of a linear form. For the fourth-order fully symmetric tensor, De Lathauwer in \cite{Lieven} proposed the Fourth-Order-Only Blind Identification (FOOBI) algorithm. 

Among those numerical algorithms, the ALS method is the most popular one since it is robust. However, the ALS has some drawbacks. For example, the convergence of ALS can be extremely slow. In addition, the ALS method for SOPD is not efficient since all three subproblems are the same equation and subproblems are now nonlinear in factor matrices corresponding to the symmetry. There are very few numerical methods for finding SOPD. Schultz \cite{Schultz} numerically solves SOPD using the best symmetric rank-1 approximation of a symmetric tensor through the maximum of the associated homogeneous form over the unit sphere. In this paper, we study the SOPD for the third-order partially symmetric tensors and the fourth-order fully symmetric tensors and propose a new method called Partial Column-wise Least-squares (PCLS) to solve the SOPD. It obviates the three nonlinear least-squares subproblems through some optimized matricizations and performing a root finding technique for polynomials in finding factor matrices.

\subsection{Preliminaries}

We denote the scalars in $\mathbb{R}$ with lower-case letters $(a,b,\ldots)$ and the vectors with bold lower-case letters $(\bf{a},\bf{b},\ldots)$.  The matrices are written as bold upper-case letters $(\bf{A}, \bf{B},\ldots)$ and the symbols for tensors are calligraphic letters $(\mathcal{A},\mathcal{B},\ldots)$. The subscripts represent the following scalars:  $\mathcal{(A)}_{ijk}=a_{ijk}$, $(\bold{A})_{ij}=a_{ij}$, $(\bold{a})_i=a_i$ and the $r$-th column of a matrix $\bold{A}$ is $\bold{a_r}$. The matrix sequence is $\{\bold{A}^k\}$. 

Here we describe several necessary definitions.

\begin{definition}[Mode-$n$ matricization]
Matricization is the process of reordering the elements of an $N$th order tensor into a matrix. The mode-$n$ matricization of a tensor $\mathcal{T} \in \mathbb{R}^{I_1 \times I_2 \times \cdots \times I_N}$ is denoted by $\bold{T_{(n)}}$ and arranges the mode-$n$ fibers to be the columns of the resulting matrix.  The mode-$n$ fiber, $\bold{t_{i_1\cdots i_{n-1} : i_{n+1} \cdots i_N}}$, is a vector obtained by fixing every index with the exception of the $n$th index.
\end{definition}

If we use a map to express such matricization process for any $N$th order tensor $\mathcal{T} \in \mathbb{R}^{I_1 \times I_2 \times \cdots \times I_N}$, that is, the tensor element $(i_1, i_2, \dots, i_N)$ maps to matrix element $(i_n, j)$, then there is a formula to calculate $j$: 
$$j=1+\sum_{\substack{k=1\\ k\neq n}}^{N}(i_k -1)J_k \quad \text{with} \quad J_k =\prod_{\substack{m=1\\ m\neq n}}^{k-1}I_m.$$
So, given a third-order tensor $\mathcal{X} \in \mathbb{R}^{I \times J \times K}$, the mode-$1$, mode-$2$ and mode-$3$ matricizations of $\mathcal{X}$, respectively, are:
\begin{eqnarray}\label{matricize}
  \bold{X_{(1)}}&=&[\bold{x_{:11}}, \dots, \bold{x_{:J1}}, \bold{x_{:12}}\dots,\bold{x_{:J2}}, \dots, \bold{x_{:1K}}, \dots, \bold{x_{:JK}}] , \nonumber\\
  \bold{X_{(2)}}&=&[\bold{x_{1:1}}, \dots, \bold{x_{I:1}}, \bold{x_{1:2}}\dots,\bold{x_{I:2}}, \dots, \bold{x_{1:K}}, \dots, \bold{x_{I:K}}] ,\\
   \bold{X_{(3)}}&=&[\bold{x_{11:}}, \dots, \bold{x_{I1:}}, \bold{x_{12:}}\dots,\bold{x_{I2:}}, \dots, \bold{x_{1J:}}, \dots, \bold{x_{IJ:}}]. \nonumber
\end{eqnarray}
 
 
 \begin{definition}[square matricization]
 For a fourth-order tensor $\mathcal{T} \in \mathbb{R}^{I \times J \times K \times L}$, the square matricization is denoted by $mat(\mathcal{T}) \in \mathbb{R}^{IK \times JL}$ and is defined as
 \begin{eqnarray}\label{eq:sqmat}
 \bold{T} = mat(\mathcal{T}) \Leftrightarrow (\bold{T})_{(i-1)K +k, (j-1)L+l} = \mathcal{T}_{ijkl}.
 \end{eqnarray}
\end{definition}
See the paper \cite{teninversion} for the generalizations of square matricication in terms of tensor blocks.

{\begin{definition}[unvec]
Given a vector $\bold{v} \in \mathbb{R}^{I^2}$, $unvec(\bold{v}) = \bold{W}$ is a square matrix of size $I \times I$ obtained from matricizing $\bold{v}$ via through its column vectors $\bold{w}_j \in \mathbb{R}^I$, $j= 1, 2, \dots, I$; i.e. 
$$\bold{w}_{ij} =v((j-1)\cdot I+i),~~i=1,2,\ldots,I$$
and  
$$unvec(\bold{v}) = \begin{bmatrix} \bold{w}_1 & \bold{w}_2 & \dots & \bold{w}_I \end{bmatrix}.$$
 \end{definition}


\section{Symmetric Outer Product Decomposition}

\begin{definition}
Let $\bold{x},\bold{y} \in \mathbb{R}^n$. The outer product of $\bold{x}$ and $\bold{y}$ is
\begin{eqnarray}
 \bold{M}=\left [
\begin{array}{cccc}
x_{1} y_{1} & x_{1} y_{2} & \cdots & x_{1} y_{n}  \\
x_{2} y_{1} &                     &               & \vdots \\
\vdots          &                     &               & \vdots\\
x_{n} y_{1} &                     &                & y_{n} y_{n} 
\end{array}
\right ].
\end{eqnarray}
\end{definition}

If $\bold{x}=\bold{y}$, then we see that $\bold{M}$ is a symmetric matrix.

\noindent
The outer product of the vectors $\bold{x}, \bold{y}, \bold{z}  \in \mathbb{R}^n$ is the following:

\begin{eqnarray}
(\bold{x} \otimes \bold{y} \otimes \bold{z})_{ijk} = x_i y_j z_k.
\end{eqnarray}

\noindent
The outer product of three vectors is a third-order rank-one tensor; the outer product of $k$ vectors is a $k$th-order rank-one tensor.
Let $\mathcal{T}=\bold{x} \otimes \bold{y} \otimes \bold{z}$, moreover, if $\bold{x}=\bold{y}=\bold{z}$, then we say $\mathcal{T}$ is a symmetric third-order rank-one tensor.
If either $\bold{x}=\bold{y}$, $\bold{x}=\bold{z}$ or $\bold{y}=\bold{z}$, then we say $\mathcal{T}$ is a partially symmetric third-order rank-one tensor.

\begin{definition}[Rank-one tensor] 
A $k$th order tensor $\mathcal{T} \in \mathbb{R}^{\mathnormal{I}_1 \times \mathnormal{I}_2 \times \cdots \times \mathnormal{I}_k}$ is called rank-one if it can be written as an outer product of $k$ vectors; i.e.
$$\mathcal{T}_{i_1 i_2 \cdots i_k} = a_{i_1}^{(1)}a_{i_2}^{(2)} \cdots a_{i_k}^{(k)}, \quad \text{for all} \; \; 1 \leq i_r \leq I_r.$$
Conveniently,  a rank-one tensor is expressed as
$$\mathcal{T} = \bold{a}^{(1)} \otimes  \bold{a}^{(2)} \otimes \cdots \otimes  \bold{a}^{(k)},$$
where $\bold{a}^{(r)} \in \mathbb{R}^{I_r}$ with $1 \leq r \leq k$. 
\end{definition}

\begin{definition}[Symmetric rank-one tensor] \label{ranksym}
A rank-one $k$th-order tensor $\mathcal{T} \in \R^{I \times I \times \cdots \times I}$ is symmetric if it can be written as an outer product of $k$ vectors; i.e.
$$\mathcal{T} = \underbrace{ \bold{a} \otimes \bold{a} \otimes \cdots \otimes \bold{a}}_{\text{$k$}}$$
where $\bold{a} \in \R^{I}$.
\end{definition} 


\begin{remark} We say a tensor is cubical if its modal dimensions are identical. Symmetric tensors are cubical.
A fully symmetric tensor is invariant under all permutations of its indices.  Let the permutation $\sigma$ be defined as $\sigma(i_1,i_2,\ldots,i_k)=i_{m(1)}i_{m(2)} \ldots i_{m(k)}$ where $m(j) \in \{1,2,\ldots,k\}$. If $\mathcal{T}$ is a symmetric tensor, then
\[    \mathcal{T}_{\sigma}=\mathcal{T}_{i_{m(1)}i_{m(2)} \ldots i_{m(k)} }  \]
for all permutation $\sigma$ on the index set $\{i_1,i_2,\ldots,i_k\}$.
\end{remark}


\begin{definition}[Partially symmetric rank-one tensor] \label{rankpsym}
A rank-one $k$th-order tensor $\mathcal{T} \in \R^{I_1 \times I_2 \times \cdots \times I_k}$ is partially symmetric if it can be written as an outer product of $k$ vectors and if
there exist modes $l$ and $m$ such that
$\bold{a}^{(l)} = \bold{a}^{(m)}$ where $1 \leq l, m \leq k$ and $l \neq m$ in
$$\mathcal{T} = \bold{a}^{(1)} \otimes \bold{a}^{(2)} \otimes  \ldots \otimes \bold{a}^{(k)}  $$
with $\bold{a}^{(r)} \in \R^{I_r}$.
\end{definition}
Above is a minimal definition for a tensor to have partial symmetry.  There can exist disjoint subindices, $S_1=\{s^1_i, i \in \{1,2,\ldots,k  \} \}, S_2=\{s^2_i, i \in \{1,2,\ldots,k \} \}, \ldots$ for which $\bold{a}^{(s^1_1)} = \bold{a}^{(s^1_2)} = \cdots = \bold{a}^{(s^1_{|S_1|})}$, $\bold{a}^{(s^2_1)} = \bold{a}^{(s^2_2)} = \cdots = \bold{a}^{(s^2_{|S_2|})}$ and etc.


\begin{remark}
If a third-order tensor $\mathcal{T}$ is partially symmetric tensor with $\bold{a}^{(1)}= \bold{a}^{(2)}$, then
\[    \mathcal{T}_{i_1 i_2 i_3}=\mathcal{T}_{i_2 i_1 i_3}.  \]
\end{remark}

A $k$th-order tensor $\mathcal{T}$ can be decomposed into as sum of outer products of vectors if there exists a positive number $R$ such that
\begin{eqnarray}
\mathcal{T}=\sum_{r=1}^{R} \underbrace{ \bold{a}_r^{(1)} \otimes \bold{a}_r^{(2)} \otimes \cdots \otimes \bold{a}_r^{(k)}}_{\text{$k$}}
\end{eqnarray}
exists. This is called the Canonical Polyadic (CP) decomposition (also known as PARAFAC and CANDECOM). This decomposition into a sum of a symmetric and/or unsymmetric outer product decompositions first appeared in the papers of Hitchcock \cite{Hitch1,Hitch2}.  The notion of tensor rank was also introduced by Hitchcock.

\begin{definition}\label{}
The rank of $\mathcal{T} \in \mathbb{R}^{I_1 \times \cdots \times I_k}$ is defined as
\[ rank(\mathcal{T}):= \min_{R} \Big \{  R \Big | \mathcal{T}=\sum_{r=1}^{R}  \bold{a}_r^{(1)} \otimes \bold{a}_r^{(2)} \otimes \cdots \otimes \bold{a}_r^{(k)}   \Big \} \]
\end{definition}

Define $\mathsf{T}^k(\mathbb{R}^{n})$ as the set of all order-$k$ dimensional $n$ cubical tensors. A set of symmetric tensors in $\mathsf{T}(\mathbb{R}^{n})$ is denoted as
$\mathsf{S}^k(\mathbb{R}^n)$.

\begin{definition}\label{}
If $\mathcal{T} \in \mathsf{S}^k(\mathbb{R}^n)$, then the rank of a symmetric  $\mathcal{T} \in \mathbb{R}^{I_1 \times \cdots \times I_k}$ is defined as
\[ rank_{\mathsf{S}}(\mathcal{T}):= \min_{S} \Big \{  S \Big | \mathcal{T}=\sum_{s=1}^{S}  \underbrace{\bold{a}_s \otimes \bold{a}_s  \otimes \cdots \otimes \bold{a}_s}_{k}   \Big \} \]
\end{definition}


\begin{lemma}\cite{Comon1}
Let $\mathcal{T} \in \mathsf{S}^k(\mathbb{R}^n)$, there exist $\bold{x}_1, \bold{x}_2, \cdots, \bold{x}_S \in \mathbb{R}^n$ linearly independent vectors such that
\[   \mathcal{T}= \sum_{i=1}^S   \underbrace{\bold{x}_i\otimes \bold{x}_i \otimes \cdots \otimes \bold{x}_i}_{k}  \]
has $rank_{\mathsf{S}}( \mathcal{T})=S$.
\end{lemma}

Note that $\mathsf{S}^k(\mathbb{R}^n) \subset \mathsf{T}^k(\mathbb{R}^{n})$.
We have that $R(k,n) \geq R_{\mathsf{S}}(k,n)$ where $R(k,n)$ be the maximally attainable rank in the space of order-$k$ dimension-$n$ cubical tensors $\mathsf{T}^k(\mathbb{R}^{n})$ and
$R_{\mathsf{S}}(k,n)$ be the maximally attainable symmetric rank in the space of symmetric tensors $\mathsf{S}^k(\mathbb{R}^n)$. In \cite{Comon1,Landsberg}, there are numerous results on symmetric rank over $\mathbb{C}$. For example in \cite{Comon1}, for all $\mathcal{T}$
\begin{itemize}
\item
$rank_{\mathsf{S}}(\mathcal{T}) \leq \binom{n+k-1}{k}$
\item
$rank(\mathcal{T}) \leq rank_{\mathsf{S}}(\mathcal{T})$
 \end{itemize}
We also refer the readers to the book by Landsberg \cite{Landsberg} on some discussions on partially symmetric tensor rank and the work of Stegeman \cite{AStegeman} on some uniqueness conditions for the minimum rank of symmetric outer product.


%


\section{Alternating Least-Squares}\label{sec:three}

Our goal is approximating a minimum sum of rank-one $k$th-order tensors from a given tensor $\mathcal{T}$. Given a $k$th-order tensor $\mathcal{T} \in \mathbb{R}^{I_1 \times I_2 \times \ldots \times I_k}$, 
find the best minimum sum of rank-one $k$th-order tensor 
\begin{eqnarray}\label{bestcpformulate}
\min_{R} \Vert \mathcal{T} - \mathcal{\widetilde{T}}   \Vert_F^2
\end{eqnarray}
where  $\mathcal{\widetilde{T}}=\displaystyle{\sum_{r=1}^R \bold{a}_r^{(1)} \otimes \bold{a}_r^{(2)}} \otimes \cdots \otimes \bold{a}_r^{(k)}$.


ALS is a numerical method for approximating the canonical decomposition of a given tensor.
For simplicity, we describe ALS for third-order tensors.  The ALS problem for third order tensor is the following
\begin{eqnarray*}
\displaystyle\mathop{\mathrm{min}}_{\bold{A},\bold{B},\bold{C}} \quad \left\Vert\mathcal{T}-\sum_{r=1}^{R}\bold{a}_r \otimes \bold{b}_r \otimes \bold{c}_r\right\Vert_F^2
\end{eqnarray*}
where $\mathcal{T} \in \mathbb{R}^{I \times J \times K}$. Define the factor matrices $\bold{A}$, $\bold{B}$ and $\bold{C}$ as the concatenation of the vectors $\bold{a}_r$, $\bold{b}_r$ and  $\bold{c}_r$, respectively; i.e.,
$\bold{A}=[\bold{a}_1~\bold{a}_2~ \ldots \bold{a}_R] \in \mathbb{R}^{I \times R}$, $\bold{B}=[\bold{b}_1~\bold{b}_2~ \ldots \bold{b}_R] \in \mathbb{R}^{J \times R}$ and $\bold{C}=[\bold{c}_1~\bold{c}_2~ \ldots \bold{c}_R]  \in \mathbb{R}^{K \times R}$.


Matricizing the equation
\[   \mathcal{T} =\sum_{r=1}^{R}\bold{a}_r \otimes \bold{b}_r \otimes \bold{c}_r   \]
on both sides, we obtain three equivalent matrix equations:
\begin{eqnarray*}
&&\bold{T_{(1)}} = \bold{A}(\bold{C}\odot\bold{B})^{\text{T}}, \\
&&\bold{T_{(2)}} = \bold{B}(\bold{C}\odot\bold{A})^{\text{T}}, \\
&&\bold{T_{(3)}} = \bold{C}(\bold{B}\odot\bold{A})^{\text{T}}.
\end{eqnarray*}
where  $\bold{T_{(1)}}^{I\times JK}$, $\bold{T_{(2)}}^{J\times IK}$ and $\bold{T_{(3)}}^{K\times IJ}$ are the mode-1, mode-2 and mode-3 matricizations of tensor $\mathcal{T}$. The symbol $\odot$ denotes the Khatri-Rao product \cite{RM71}.  Given matrices $\bold{A} \in \mathbb{R}^{\mathnormal{I} \times \mathnormal{R}}$ and $\bold{B} \in \mathbb{R}^{\mathnormal{J} \times \mathnormal{R}}$, the Khatri-Rao product of $\bold{A}$ and $\bold{B}$ is the ``matching columnwise" Kronecker product; i.e., 
\begin{eqnarray*}
\bold{A} \odot \bold{B}= [\bold{a_1} \otimes \bold{b_1}~~ \bold{a_2} \otimes \bold{b_2}~~ \ldots] \in \mathbb{R}^{IJ \times K}.
\end{eqnarray*}
By fixing two factor matrices but one at each minimization, three coupled linear least-squares subproblems are then formulated to find each factor matrices:
\begin{eqnarray}\label{als}
\bold{A}^{k+1}&=&\displaystyle\mathop{\mathrm{argmin}}_{\widehat{\bold{A}}\in \mathbb{R}^{I\times R}} \left\Vert\bold{T_{(1)}}^{I\times JK}-\widehat{\bold{A}}  (\bold{C}^{k}\odot\bold{B}^{k})^{\text{T}}\right\Vert_F^2, \notag \\
\bold{B}^{k+1}&=&\displaystyle\mathop{\mathrm{argmin}}_{\widehat{\bold{B}}\in \mathbb{R}^{J\times R}}\left\Vert\bold{T_{(2)}}^{J\times IK}-\widehat{\bold{B}}(\bold{C}^{k}\odot\bold{A}^{k+1})^{\text{T}}\right\Vert_F^2, \\
\bold{C}^{k+1}&=&\displaystyle\mathop{\mathrm{argmin}}_{\widehat{\bold{C}}\in \mathbb{R}^{K\times R}}\left\Vert\bold{T_{(3)}}^{K\times IJ}-\widehat{\bold{C}}(\bold{B}^{k+1}\odot\bold{A}^{k+1})^{\text{T}}\right\Vert_F^2.  \notag
\end{eqnarray}
where $\bold{T_{(1)}}$, $\bold{T_{(2)}}$ and $\bold{T_{(3)}}$ are the standard tensor flattennings described in (\ref{matricize}).
To start the iteration, the factor matrices are initialized with $\bold{A}^0$, $\bold{B}^0$, $\bold{C}^0$. ALS fixes $\bold{B}$ and $\bold{C}$ to solve for $\bold{A}$, then it fixes $\bold{A}$ and $\bold{C}$ to solve for $\bold{B}$.  And then ALS finally fixes $\bold{A}$ and $\bold{B}$ to solve for $\bold{C}$. This Gauss-Seidel sweeping process continues iteratively until some convergence criterion is satisfied. Thus the original nonlinear optimization problem can be solved with three linear least squares problems.

ALS is a very simple method that has been applied across engineering and science disciplines. However, ALS has some disadvantages. For non-degenerate problems, convergence may require a high number of iterations (see Figure \ref{swamp1}) which can be attributed to the non-uniqueness in the solutions of the subproblems, collinearity of the columns in the factor matrices and initialization of the factor matrices; see e.g. \cite{Comon5,Rajih,Paatero}.  The long curve in the residual plot is also an indication of a degeneracy problem. 

The ALS algorithm can be applied to find symmetric and partially symmetric outer product decomposition for third order tensor by setting $\bold{A}=\bold{B}=\bold{C}$ and  $\bold{A} = \bold{B}$ or $\bold{A} = \bold{C}$, respectively, in (\ref{als}). But the swamps are prevalent in these cases and the factor matrices obtained often do not reflect the symmetry of the tensor. In addition, when ALS is applied to symmetric tensors, the least-squares subproblems are highly ill-conditioned which lead to non-unique solutions in all three directions. As ALS cycles through the iterations, these subproblems pull together to drive the outputs away from the true solutions. The regularization methods \cite{RALS1,RALS2} does not drastically alleviate this type of swamps.

\begin{figure}[h]\label{swamp1}
\centering
\includegraphics[width=6.3cm]{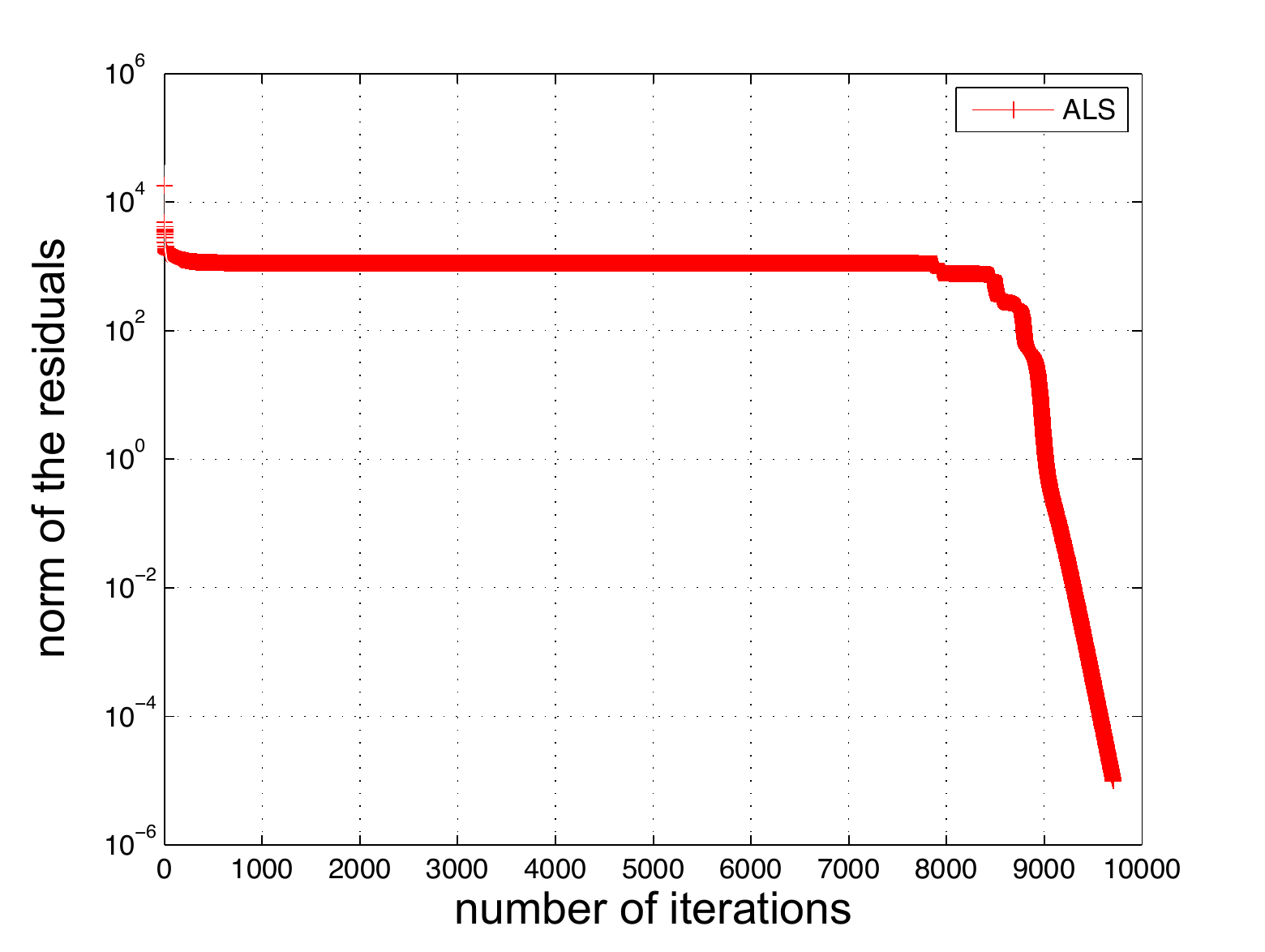}
\caption{The long flat curve (\emph{swamp}) in the ALS method. The error stays at $10^3$ during the first 8000 iterations. }
\end{figure}

%


\noindent
Here are the problem formulations: given an order-$k$th tensor $\mathcal{T} \in \mathbb{R}^{I_1 \times I_2 \times \ldots \times I_k}$, 

\begin{enumerate}
\item[(1)]
find the best minimum sum of rank-one \emph{symmetric} tensor 

\begin{eqnarray*}
\min_{} \Vert \mathcal{T} - \mathcal{\widetilde{T}}   \Vert_F^2~~~~~~\mbox{(Problem~1)}
\end{eqnarray*}

where $\mathcal{\widetilde{T}}= \displaystyle{\sum_{r=1}^{R} \bold{a}_r \otimes \bold{a}_r \otimes \cdots \otimes \bold{a}_r }$

\item[(2)]
find the best minimum sum of rank-one \emph{partially symmetric} tensor

\begin{eqnarray*}
\min_{} \Vert \mathcal{T} - \mathcal{\widetilde{T}} \Vert_F^2~~~~~~\mbox{(Problem~2)}
\end{eqnarray*}
where $\mathcal{\widetilde{T}}=\displaystyle{\sum_{r=1}^{R}\bold{a}_r^{(1)} \otimes \bold{a}_r^{(2)} \otimes  \ldots \otimes \bold{a}_r^{(k)}} $ for some modes $\bold{a}_r^{(j)}=\bold{a}_r^{(l)}$ where $1 \leq j, l \leq k$ and $j \neq l$.
\end{enumerate}
We refer to these decomposition symmetric outer product decompositions (SOPD).\\

For the sake of clarity of the exposition, we describe the decomposition methods for third-order and forth-order tensors with partial and full symmetries. We also include some discussions on how these methods can be extended  to the general case.


\subsection{SOPD for Third-order Partially Symmetric Tensor}

Given a third-order tensor  $\mathcal{T} \in \R^{I \times I \times K}$ with $t_{ijk} = t_{jik}$, Problem 2 becomes

\begin{eqnarray}\label{mainprob1}
\displaystyle\mathop{\mathrm{min}}_{\bold{A},\bold{C}} \quad \left\Vert\mathcal{T}-\sum_{r=1}^{R_{ps}}\bold{a}_r \circ \bold{a}_r \circ \bold{c}_r\right\Vert_F^2,
\end{eqnarray}
with $R_{ps}$ summands of rank-one partial symmetric tensors and $\widehat{\mathcal{T}}=\sum_{r=1}^{R_{ps}}\bold{a}_r \circ \bold{a}_r \circ \bold{c}_r$ . The unknown vectors are arranged into two factor matrices $\bold{A} = [ \bold{a}_1 ~ \bold{a}_2 ~ \cdots ~ \bold{a}_{R_{ps}} ]$ and $\bold{C} = [\bold{c}_1 ~ \bold{c}_2 ~ \cdots ~ \bold{c}_{R_{ps}}]$ in this case. Matricization of $\widehat{\mathcal{T}}$ leads to

\begin{eqnarray*}
\bold{\widehat{T}_{(3)}} &=& \bold{C} (\bold{A} \odot \bold{A})^{\text{T}},
\end{eqnarray*}
where $\bold{\widehat{T}_{(3)}} \in \R^{K \times I^2}$ is the mode-3 matricization of tensor $\widehat{\mathcal{T}}$. Thus (\ref{mainprob1}) becomes
\begin{eqnarray}\label{mainprob2}
\displaystyle\mathop{\mathrm{min}}_{\bold{A},\bold{C}}\quad \left\Vert \bold{T_{(3)}} - \bold{C}(\bold{A} \odot \bold{A})^{\text{T}} \right\Vert_F^2.
\end{eqnarray}
If we apply the ALS method, the problem reduces to the following subproblems:
\begin{eqnarray}
\bold{A}^{k+1} &=& \displaystyle\mathop{\mathrm{argmin}}_{\widehat{\bold{A}} \in \R^{I \times R_{ps}}} \left\Vert \bold{T_{(3)}} - \bold{C}^k(\widehat{\bold{A}} \odot \widehat{\bold{A}})^{\text{T}} \right\Vert_F^2, \label{main-eq1} \\
\bold{C}^{k+1} &=& \displaystyle\mathop{\mathrm{argmin}}_{\widehat{\bold{C}} \in \R^{K \times R_{ps}}} \left\Vert \bold{T_{(3)}} - \bold{\widehat{C}}(\bold{A}^{k+1} \odot \bold{A}^{k+1})^{\text{T}} \right\Vert_F^2. \label{main-eq2}
\end{eqnarray}
Observe that (\ref{main-eq2}) is a linear least-squares subproblem, but \eqref{main-eq1} is a nonlinear least-squares subproblem. Directly applying the ALS method to equations \eqref{main-eq1} and \eqref{main-eq2} does not work; it often leads a wrong  solution; i.e., the factor matrices do not satisfy tensor symmetries, and it takes a high number of iterations (swamps) for it to converge.

To obviate this problem, we focus on (\ref{main-eq1}) and find an alternative method to solve for the factor matrix $\bold{A}$. 
Recall that $\bold{T_{(3)}} = \bold{C}^k(\widehat{\bold{A}} \odot \widehat{\bold{A}})^{\text{T}}$ can be solved for 
$\widehat{\bold{A}} \odot \widehat{\bold{A}}$; i.e. 
\begin{eqnarray} \label{mpinv}
\widehat{\bold{A}} \odot \widehat{\bold{A}} = ((\bold{C}^k)^{\dagger}\bold{T_{(3)}})^{\text{T}}
\end{eqnarray}
where  $(\cdot)^{\dagger}$ denotes the Moore-Penrose pseudoinverse. Equivalently, \eqref{mpinv} can be written as
\begin{eqnarray}
 \bold{\widehat{a}}_r \otimes  \bold{\widehat{a}}_r = ((\bold{C}^k)^{\dagger}\bold{T_{(3)}})^{\text{T}}(:,r) \Leftrightarrow \bold{\widehat{a}}_r \cdot \bold{\widehat{a}}_r^{\text{T}} = unvec\left(((\bold{C}^k)^{\dagger}\bold{T_{(3)}})^{\text{T}}(:,r)\right)  \label{target}
\end{eqnarray}
where  $r=1, 2, \dots, R_{ps}$, $\bold{\widehat{a}}_r$ is the $r$th column of matrix $\bold{\widehat{A}}$ and $unvec\left(((\bold{C}^k)^{\dagger}\bold{T_{(3)}})^{\text{T}}(:,r)\right)$ is a matrix $(I \times I)$ obtained from the vector $((\bold{C}^k)^{\dagger}\bold{T_{(3)}})^{\text{T}}(:,r)$ via column vector stacking of size $I$. With \eqref{target}, we can obtain  $\widehat{\bold{A}}$ by calculating each of its column $\bold{\widehat{a}}_r$ at a time.

Let $\bold{x} \in \R^{I} = [x_1~x_2~\cdots~x_I]^{\text{T}}$ denote the unknown vector $\bold{\widehat{a}}_r$ and $\bold{Y}=unvec\left(((\bold{C}^k)^{\dagger}\bold{T_{(3)}})^{\text{T}}(:,r)\right) \in \R^{I \times I}$. Then (\ref{target}) becomes
\begin{eqnarray*}
\begin{bmatrix} x_1^2 & x_1x_2 & \cdots & x_1x_I \\ x_1x_2 & x_2^2 & & \\ \vdots & & \ddots& \\ x_1x_I & & & x_I^2 \end{bmatrix} = \bold{Y}.
\end{eqnarray*}
Notice that the unknown $x_1$ is only involved in the first column and first row, so we only take the first column and first row elements of $\bold{Y}$. Thus, the least-squares formulation for these elements is
\begin{eqnarray}\label{poly}
x^*_1=\mbox{arg}\min_{x_1}~ (y_{11}-x_1^2)^2 + \sum_{i=2}^I\left[(y_{i1}-x_{i}x_1)^2+(y_{1i}-x_ix_1)^2\right].
\end{eqnarray}
This cost function in (\ref{poly}) is a fourth-order polynomial in one variable $x_1$. Thus each component $x_i$ can be solved in the same manner of minimizing a fourth-order polynomial. 

Here are the two subproblems with two initial factor matrices $\bold{A}^0$ and $\bold{C}^0$ :
\begin{eqnarray}
\bold{a}_r^{k+1} &=&  \displaystyle\mathop{\mathrm{argmin}}_{\widehat{a}_r \in \R^{I}} \left\Vert unvec\left(((\bold{C}^k)^{\dagger}\bold{T_{(3)}})^{\text{T}}(:,r)\right)-\bold{\widehat{a}}_r \cdot \bold{\widehat{a}}_r^{\text{T}}  \right\Vert_F^2, \notag\\
 && r=1,\dots, R_{ps},\label{updatea}\\
\bold{C}^{k+1} &=& \displaystyle\mathop{\mathrm{argmin}}_{\widehat{\bold{C}} \in \R^{K \times R_{ps}}} \left\Vert \bold{T_{(3)}} - \bold{\widehat{C}}(\bold{A}^{k+1} \odot \bold{A}^{k+1})^{\text{T}} \right\Vert_F^2 \label{updatec}
\end{eqnarray}
for approximating $\bold{A}$ and $\bold{C}$. We call this method the iterative Partial Column-wise Least-Squares (PCLS). 

Starting from the initial guesses, the first subproblem is solved for each column $\bold{a}_r$ of $\bold{A}$ while $\bold{C}$ is fixed. Then in the second subproblem, we fixed $\bold{A}$ to solve for $\bold{C}$. This process continues iteratively until some convergence criterion is satisfied. 

The advantage of PCLS over ALS  is that it directly computes two factor matrices. If the ALS method is applied to this problem, then one has to update three factor matrices even though there are only two distinct factors in each iteration. In addition, a very high number of iterations is required for this ALS problem to converge and it also not guaranteed that the solution satisfies the symmetries. The ALS method solves three linear least squares problems in each iteration, while PCLS  solves one least squares and $R_{ps}$ quartic polynomials solve in one iteration. 

The operational cost of running PCLS on a third-order tensor is much less than ALS since it requires only one linear least-squares with has a computational complexity of $\mathcal{O}(n^3)$ (via SVD, QR or Cholesky factorization) and finding roots of a quartic polynomial as opposed to the number of operations for three linear least-squares.
 
%




\subsection{SOPD for Fourth-order Partially Symmetric Tensors}
We can apply PCLS on the fourth-order partial symmetric tensor. We consider two cases:

\begin{enumerate}
\item[Case 1:] Let us consider the fourth-order partially symmetric tensor $\mathcal{X} \in \mathbb{R}^{I\times J\times I\times J}$ with $x_{ijkl} = x_{kjil}$ and $x_{ijkl}=x_{ilkj}$. The problem is to find factor matrices $\bold{A}$ and $\bold{B}$ through the following minimization
\begin{eqnarray}\label{fourthprob}
\min_{\bold{A}, \bold{B}} \left\Vert \mathcal{X} - \sum_{r=1}^R \bold{a}_r\circ \bold{b}_r \circ \bold{a}_r \circ \bold{b}_r \right\Vert_F^2,
\end{eqnarray}
where $\bold{A} = [\bold{a}_1~\bold{a}_2~\cdots ~\bold{a}_R]$ and $\bold{B} = [\bold{b}_1~\bold{b}_2~\cdots~\bold{b}_R]$.

By using the square matricization, we obtain
\begin{eqnarray}\label{fourtheq2}
mat(\mathcal{X}) = (\bold{A} \odot \bold{A}) (\bold{B} \odot \bold{B})^T.
\end{eqnarray}

To solve the equation (\ref{fourtheq2}) for $\bold{A}$ and $\bold{B}$, we apply the least squares method on
\begin{eqnarray}\label{decomp2}
\bold{a}_r \otimes \bold{a}_r = mat(\mathcal{X})((\bold{B} \odot \bold{B})^T)^{\dagger}(:,r), \quad r=1,\dots R \\
\bold{b}_r \otimes \bold{b}_r = mat(\mathcal{X})^T((\bold{A} \odot \bold{A})^T)^{\dagger}(:,r), \quad r=1,\dots R. 
\end{eqnarray}
iteratively. The two equations above can be solved by the same method in Section $3.1$. Again, we only need to solve the global minima of fourth-order polynomials. 

\item[Case 2:] Let us consider the fourth-order partially symmetric tensor $\mathcal{X} \in \mathbb{R}^{I \times J \times I \times K}$ with $x_{ijkl} = x_{kjil}$. This means our tensor is partial symmetric in mode one and mode three. The problem is to find factor matrices $\bold{A}$, $\bold{B}$ and $\bold{C}$ via
\begin{eqnarray}\label{fourthprob2}
\min_{\bold{A}, \bold{B}, \bold{C}} \left\Vert \mathcal{X} - \sum_{r=1}^R \bold{a}_r\circ \bold{b}_r \circ \bold{a}_r \circ \bold{c}_r \right\Vert_F^2,
\end{eqnarray}
where $\bold{A} = [\bold{a}_1~\bold{a}_2~\cdots ~\bold{a}_R]$, $\bold{B} = [\bold{b}_1~\bold{b}_2~\cdots~\bold{b}_R]$ and $\bold{C} = [\bold{c}_1~\bold{c}_2~\cdots~\bold{c}_R]$.

So by using the standard matricization and square matricization, we can have the following three equations,
\begin{eqnarray}
mat(\mathcal{X}) &=& (\bold{A} \odot \bold{A}) (\bold{B} \odot \bold{C})^{\text{T}}, \label{eq:eq1}\\
\bold{X}_{(2)} &=& \bold{B} (\bold{C} \odot \bold{A} \odot \bold{A})^{\text{T}}, \label{eq:eq2}\\
\bold{X}_{(4)} &=& \bold{C} (\bold{B} \odot \bold{A} \odot \bold{A})^{\text{T}}. \label{eq:eq3}
\end{eqnarray}
Therefore, given initial guesses $\{ \bold{A}^0, \bold{B}^0, \bold{C}^0 \}$, \eqref{eq:eq1} can be solved to obtain the update of $\bold{A}$ through the method in Section $3.1$ and equations \eqref{eq:eq2} and \eqref{eq:eq3} are solved through the least-squares to update $\bold{B}$ and $\bold{C}$ iteratively.
\end{enumerate}


To solve for the SOPD for given a higher-order partial symmetric tensor, general matricizations must be applied to the tensor. See the paper \cite{teninversion} on how tensor blocks provide matricizations which are then equal to Kathri-Rao products of factor matrices. These matricized equations inherently divide into subproblems which can be solved using least-squares or variants of PCLS. 

%
%
%

\subsection{SOPD for Fourth-order Fully Symmetric Outer Product Decomposition}


Given a fourth-order fully symmetric tensor $\mathcal{T} \in \mathbb{R}^{I \times I \times I \times I}$ with $t_{ijkl} = t_{\sigma(ijkl)}$ for any permutation $\sigma$ on the index set $\{ijkl\}$. We want to a find factor matrix $\bold{A} \in \mathbb{R}^{I \times R_{s}}$ such that
\begin{eqnarray}\label{eq:probfour}
\min_{\bold{A}} \left\Vert \mathcal{T} - \sum_{r=1}^{R_s} \bold{a}_r \circ \bold{a}_r \circ \bold{a}_r \circ \bold{a}_r \right\Vert_F^2,
\end{eqnarray}
where $\bold{A} = [\bold{a}_1~\bold{a}_2~\cdots~\bold{a}_{R_s}]$. 

By using the square matricization (\ref{eq:sqmat}), we have

\begin{eqnarray}\label{prob2}
\mathcal{T} &=& \sum_{r=1}^{R_s} \bold{a}_r \circ \bold{a}_r \circ \bold{a}_r \circ \bold{a}_r \notag \\
&\Downarrow& \notag \\
\bold{T} &=& (\bold{A} \odot \bold{A}) (\bold{A} \odot \bold{A})^{\text{T}}. 
\end{eqnarray}
Since $\mathcal{T}$ is symmetric, then $\bold{T}$ is a symmetric matrix. Then it follows that there exists a matrix $\bold{E}$ such that 
\begin{eqnarray}\label{eq:eigen}
\bold{T} = \bold{E} \bold{E}^{\text{T}}.
\end{eqnarray}
Comparing the equations \eqref{prob2} and \eqref{eq:eigen}, we know that there exists an orthogonal matrix $\bold{Q}$ such that 
\begin{eqnarray}\label{eq:sopdrela}
\bold{E} = (\bold{A} \odot \bold{A})\bold{Q},
\end{eqnarray} 
where $\bold{Q} \in \mathbb{R}^{R_s \times R_s}$ is an orthogonal matrix. In equation \eqref{eq:sopdrela}, the unknowns are $\bold{A}$ and $\bold{Q}$ while $\bold{E}$ is known. This is the same problem in the third-order partially symmetric tensor case, 
$$\bold{T_{(3)}} = \bold{C}(\bold{A} \odot \bold{A})^{\text{T}},$$
where $\bold{A}$ and $\bold{C}$ are unknown and $\bold{T_{(3)}}$ is known. Therefore, given the the initial guess matrix $\bold{A}^0$ and any starting orthogonal matrix $\bold{Q}^0$, we can update the factor matrix by following subproblems
\begin{eqnarray}
\bold{A}^{k+1} &=& \displaystyle\mathop{\mathrm{argmin}}_{\widehat{\bold{A}} \in \R^{I \times R_{s}}} \left\Vert \bold{E} - (\widehat{\bold{A}} \odot \widehat{\bold{A}})\bold{Q}^k \right\Vert_F^2, \label{main-eq1four} \\
\bold{P} &=& \displaystyle\mathop{\mathrm{argmin}}_{\widehat{\bold{Q}} \in \R^{R_s \times R_{s}}} \left\Vert \bold{E} - (\bold{A}^{k+1} \odot \bold{A}^{k+1})\bold{\widehat{Q}} \right\Vert_F^2. \notag 
\end{eqnarray}
We take the QR factorization of $\bold{P}$ to obtain an orthogonal matrix $\bold{O}$. Let
\begin{eqnarray} 
\bold{Q}^{k+1} &=& \bold{O}.  \label{main-eq2four}
\end{eqnarray}
where $\bold{P}= \bold{O} \bold{R}$ and $\bold{R}$ is an upper triangular matrix. To solve equation \eqref{main-eq1four}, we apply the PCLS \eqref{updatea} to compute $\bold{A}$ column by column,
\begin{eqnarray}\label{eq:column}
\bold{a}_r^{k+1} =  \displaystyle\mathop{\mathrm{argmin}}_{\widehat{a}_r \in \R^{I}} \left\Vert unvec\left(\bold{E}(\bold{Q}^k)^{\dagger}(:,r)\right)-\bold{\widehat{a}}_r \cdot \bold{\widehat{a}}_r^{\text{T}}  \right\Vert_F^2, r=1,\dots, R_{s}.
\end{eqnarray}

We summarize the PCLS method for fourth-order fully symmetric tensor. Given the tensor $\mathcal{T} \in \mathbb{R}^{I \times I \times I \times I}$, we first calculate matrix $\bold{E} \in \mathbb{R}^{I^2 \times R_s}$ through $\bold{T}$, the matricization of $\mathcal{T}$. Then starting from the initial guesses, we fix $\bold{Q}$ to solve for each column $\bold{a}_r$ of $\bold{A}$, then $\bold{A}$ is fixed to compute a temporary matrix $\bold{P}$. In order to make sure the updated $\bold{Q}$ is orthogonal, we apply QR factorization on $\bold{P}$ to get an orthogonal matrix and set it to be the updated $\bold{Q}$. This process continues iteratively until some convergence criterion is satisfied.

\section{Numerical Examples}\label{sec: example}

In this section, we compare the performance of ALS against PCLS for the third-order partially symmetric tensors and the fourth-order fully symmetric tensors. From these numerical examples, PCLS outperformed the ALS method with respect to the number of iterations for convergence (swamp-free) and the CPU time. 

\subsection{Example I: third-order partially symmetric tensor}\label{example1chap5}
We generate a partially symmetric tensor $\mathcal{X} \in \R^{17\times 17\times18}$ by random data, in which $x_{ijk}=x_{jik}$. Consider the SOPD of $\mathcal{X}$ with $R_{ps}=17$. So it has two different factor matrices $\bold{A} \in \R^{17 \times 17}$ and $\bold{C} \in \R^{18 \times 17}$, and the decomposition is
$$\mathcal{X} = \sum_{r=1}^{R_{ps}} \bold{a}_r \circ \bold{a}_r \circ \bold{c}_r.$$

\begin{figure}[h]
\centering
\subfloat[good initial guess]{\label{fig: normal}\includegraphics[height=4.2cm]{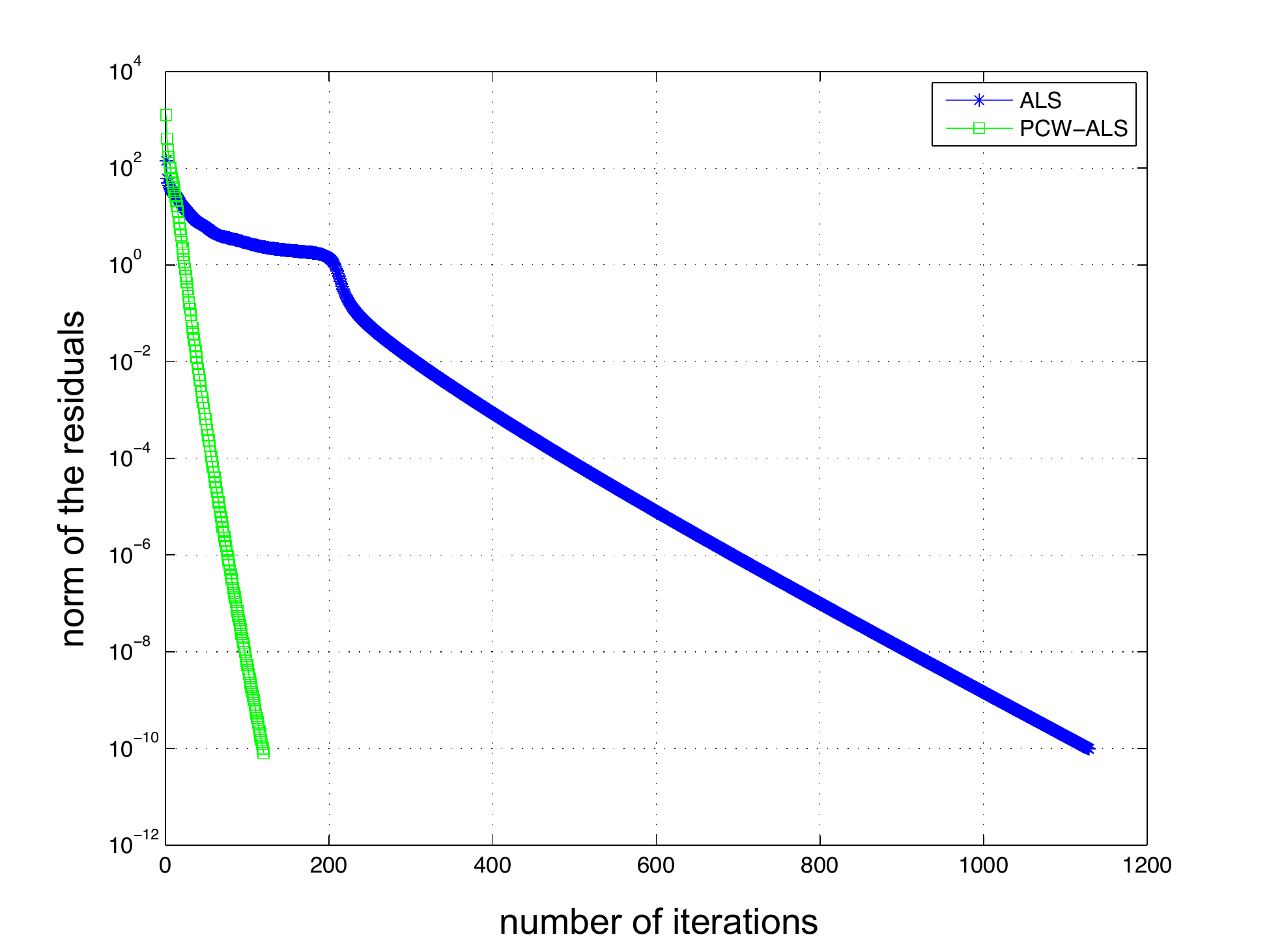}}
\qquad
\subfloat[random initial guess]{\label{fig: swamp}\includegraphics[height=4.2cm]{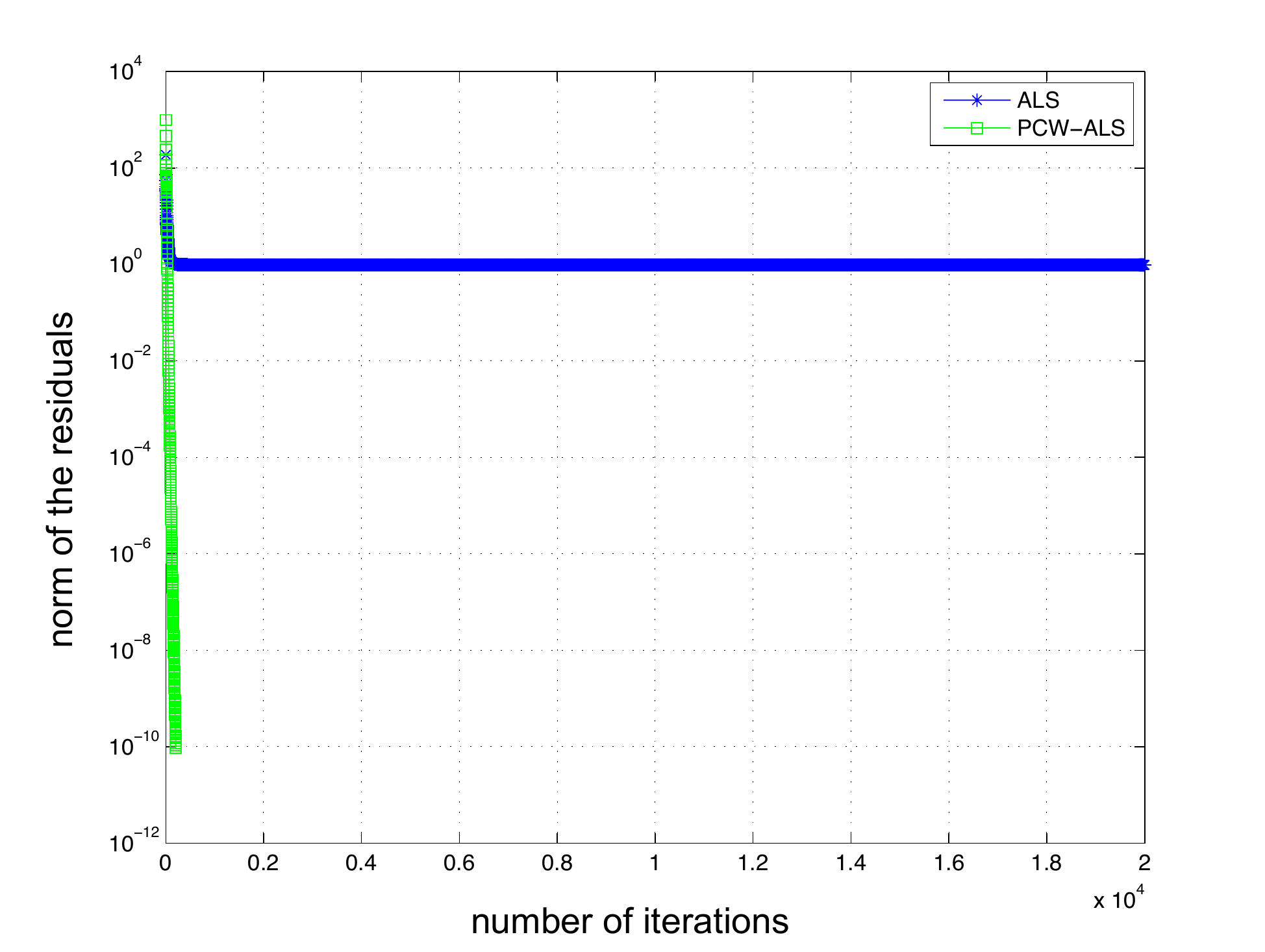}}
\caption{\bf Plots for the Example \ref{example1chap5}}
\label{fig:example1}
\end{figure}

In the two figures, the plots show the error $\left\Vert \mathcal{X}-\mathcal{X}_{est} \right\Vert_F^2$ versus the number of iterations it takes to obtain an error of $10^{-10}$, where $\mathcal{X}_{est}$ denotes the obtained tensor after every iteration. Since the ALS method needs three initial guesses, here we let $\bold{B}^0 = \bold{A}^0$ for it.

In Figure \ref{fig: normal}, the initial guesses are good. Both algorithms work well, but the PCLS method is better than the ALS algorithm. The PCLS only takes 120 iterations in comparison to that of 1129 ALS iterations. Moreover, the PCLS is faster than ALS since the CPU time of PCLS is 3.9919s while the ALS is 6.4126s. Figure \ref{fig: swamp} shows that PCLS can reduce the swamp by only taking 205 iterations to reach an error within $10^{-10}$. While the ALS has a swamp and the error stays in $10^{0}$ after 20000 iterations.

\subsection{Example II: Simulation}\label{example2chap5}
For the tensor $\mathcal{X}$ given in the Example \ref{example1chap5}, the ALS and PCLS algorithms are used to decompose it with rank $R_{ps}=17$. Both of ALS and PCLS are used on tensor $\mathcal{X}$ with 50 different random initial starters and and the average results in terms of number of iterations and CPU time are shown in the Table \ref{tab:table1chap5}. 

\begin{table}[h]
\begin{center}
\begin{tabular}{| c | c | c |}
\hline
& ALS & PCLS \\
\hline
average CPU time & 17.1546s & 6.1413s\\
\hline
average number of iterations & 3445.0 & 258.7 \\
\hline
\end{tabular}
\caption{\bf The comparison of ALS and PCLS (Mean).}
\label{tab:table1chap5}
\end{center}
\end{table}


\subsection{Example III: CPU time comparison in terms of tensor size}\label{example2bchap5}
We apply the ALS method and PCLS method on the third-order partially symmetric tensors $\mathcal{X}_1 \in \mathbb{R}^{10 \times 10 \times 10}$ with $R_{ps}=10$,  $\mathcal{X}_2 \in \mathbb{R}^{20 \times 20 \times 20}$ with $R_{ps}=20$, $\dots$, $\mathcal{X}_9 \in \mathbb{R}^{90 \times 90 \times 90}$ with $R_{ps}=90$ and compare the CPU times of both methods for the same tensor size. In order to have a fair comparison, for each tensor $\mathcal{X}_i$, we use the technique in Example \ref{example2chap5} to get the average CPU times of both methods. The following Figure \ref{fig:cpuvssize} shows that as the tensor size increases, the CPU time of ALS increases much faster than the PCLS time. 
\begin{figure}[h]
\centering
\includegraphics[height=6cm]{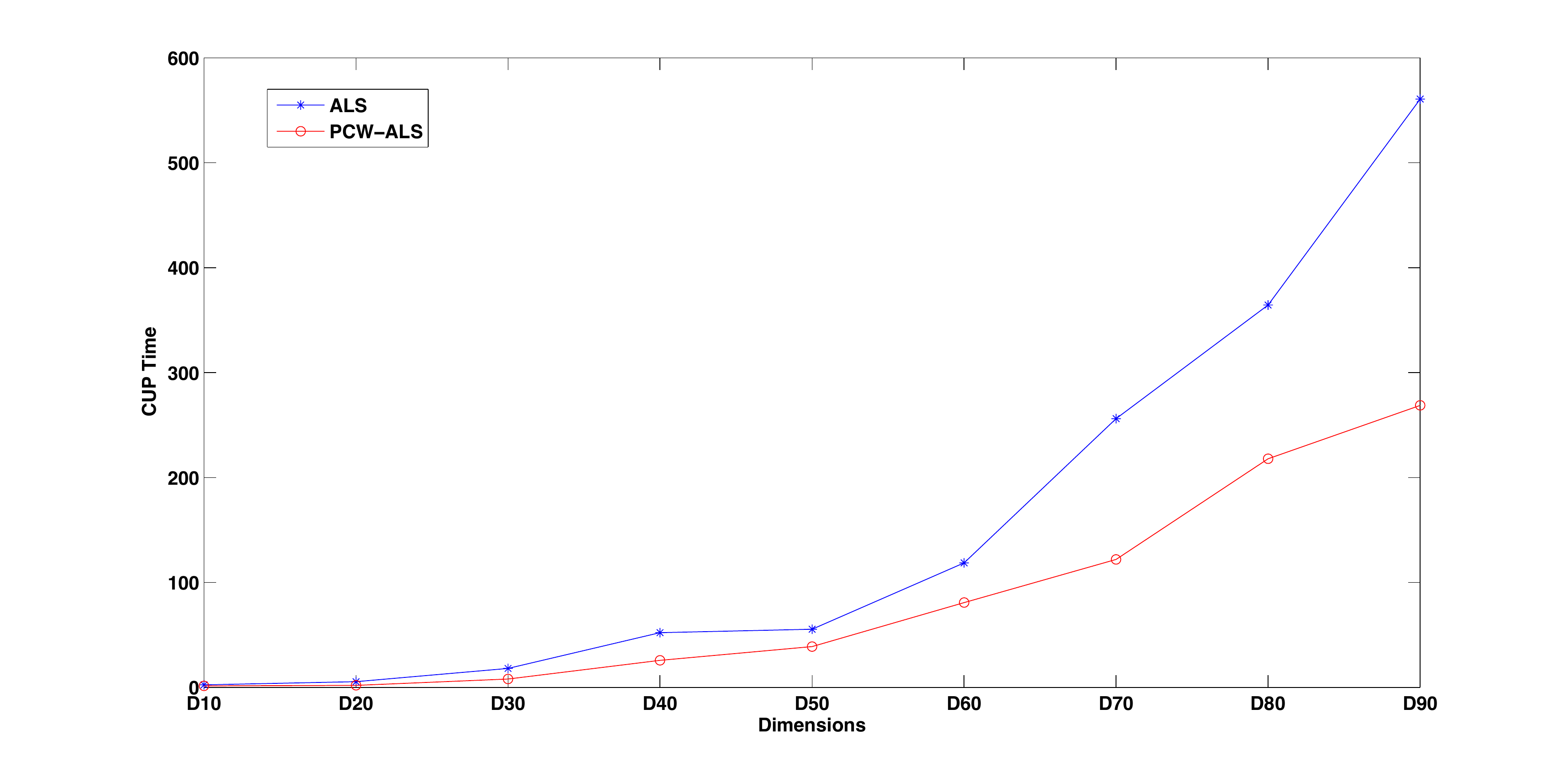}
\caption{\bf Plots for the Example \ref{example2bchap5}}
\label{fig:cpuvssize}
\end{figure}

Now we show the examples for the fourth-order fully symmetric tensor. 

\subsection{Example IV }\label{examplefourth1}
 Given fully symmetric fourth-order tensor $\mathcal{X} \in \mathbb{R}^{10 \times 10 \times 10 \times 10}$ with $R=10$, we give the initial guess $\bold{A}^0$, the ALS method and PCLS method are applied to solve the SOPD for this fourth-order tensor. The following Figure \ref{fig:figure2} shows that the swamp happens in the ALS method while the PCLS converges very fast. 
\begin{figure}[h]
\centering
\includegraphics[width=6.5cm]{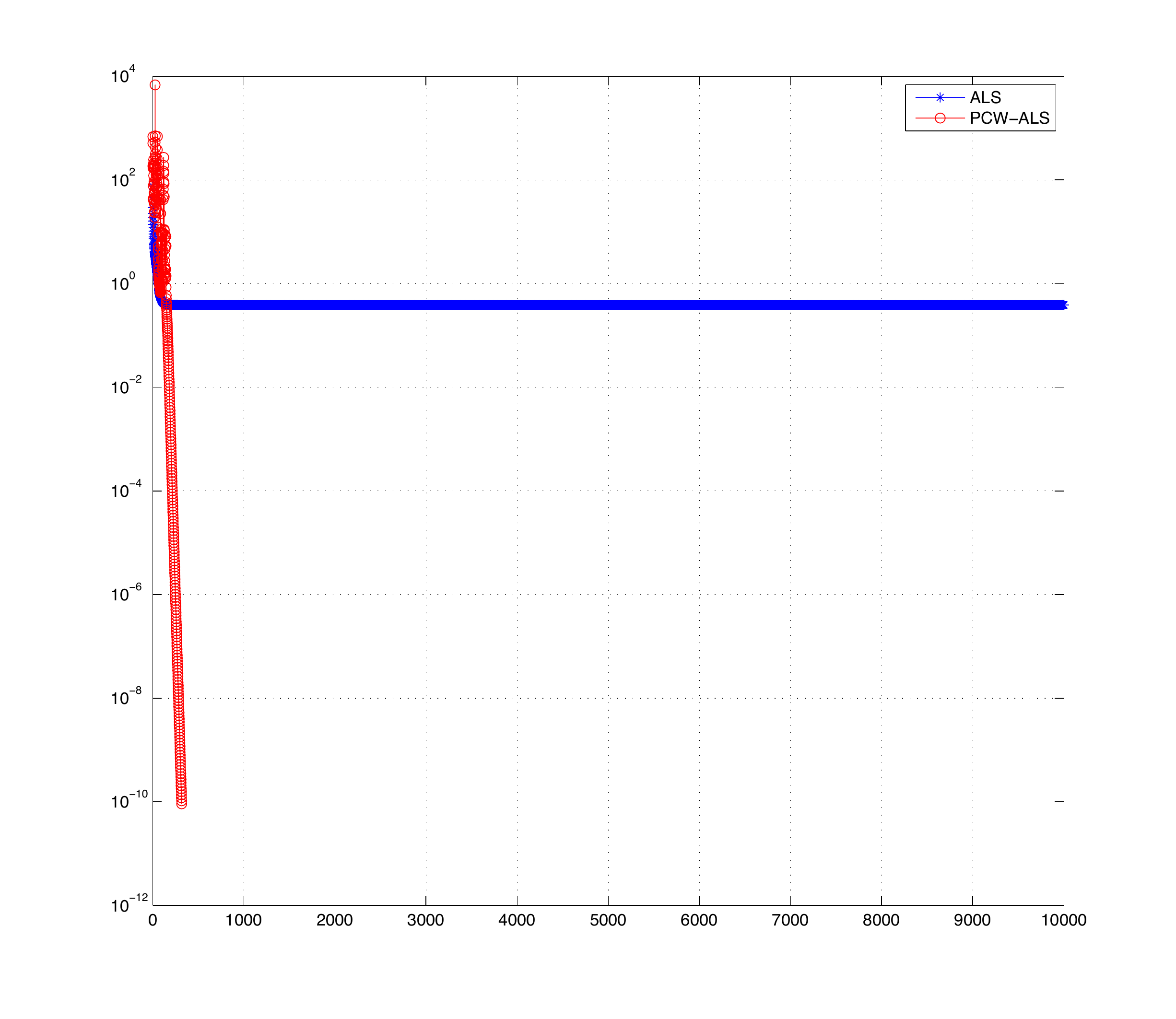}
\caption{\bf Plot for the Example \ref{examplefourth1} }
\label{fig:figure2}
\end{figure}

\subsection{Example V} \label{examplefourth2}
Given fully symmetric fourth-order tensor $\mathcal{X} \in \mathbb{R}^{15 \times 15 \times 15 \times 15}$ with $R=10$, we give the initial guess $\bold{A}^0$, the ALS method and PCLS method are applied to solve the SOPD for this fourth-order tensor. Figure \ref{fig:figure3} shows that both method works well. But the PCLS is faster than the ALS method. The CPU time of the ALS method is 27.2149s while the PCLS method is 4.2763s. 
\begin{figure}[h]
\centering
\includegraphics[width=6.5cm]{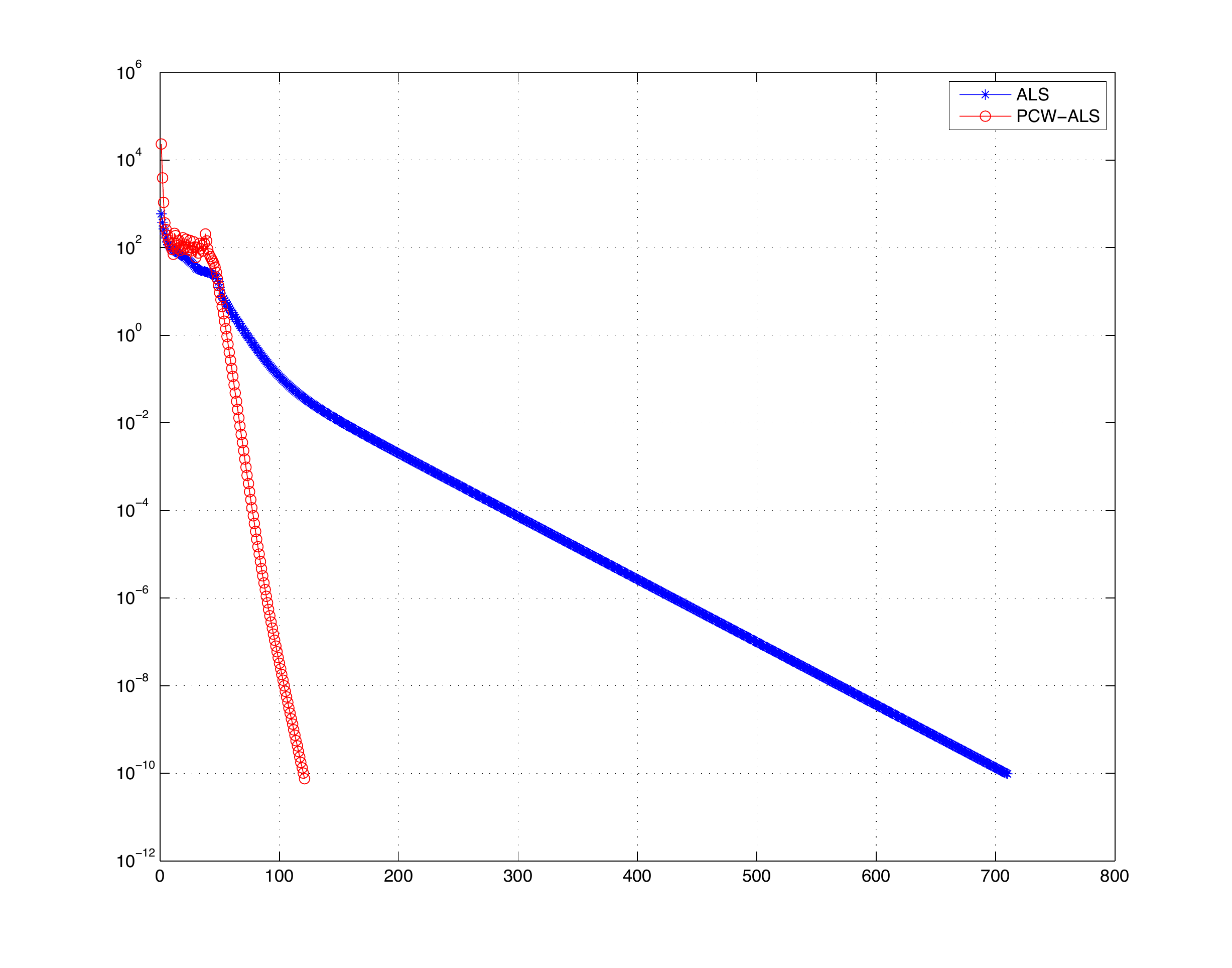}
\caption{\bf Plot for the Example \ref{examplefourth2}}
\label{fig:figure3}
\end{figure}

\section{Conclusion}\label{sec:con}

We presented the iterative algorithm PCLS for the SOPD of third-order partially symmetric tensors and  fourth-order fully symmetric tensors. The third-order partially symmetric tensor has the same factor matrix in terms of the symmetric modes, the PCLS avoided two least-squares problems for factor matrices in each iteration by solving for the roots of a quartic polynomials which updates column vectors at a time. For the fourth-order fully symmetric tensor, we reformulate the problem by using the square matricization in order to apply PCLS.  We also provided several numerical examples to compare the performance of  PCLS to ALS for the SOPD. In these examples, PCLS removes the swamps that are visible with the ALS method. 

\section*{Acknowledgements}
C.N. and N.L. were both in part supported by the U.S. National Science Foundation DMS-0915100. 

\end{document}